\documentclass[reqno, 11pt]{amsart}
\usepackage[margin=2.5cm]{geometry}
\usepackage{amscd}
\usepackage{latexsym}
\usepackage{amsmath}
\usepackage{amsthm}
\usepackage{amsfonts}
\usepackage{tikz}
\usetikzlibrary{shapes,arrows}
\numberwithin{equation}{section}

\newcommand{\vp}{\varphi}
\newcommand{\T}{\partial\mathbb{D}}

\newcommand{\ds}{\displaystyle}
\newcommand{\ol}{\overline}
\newcommand{\supp}{\text{\rm{supp}}}

\newcommand{\be}{\begin{equation}}
\newcommand{\ee}{\end{equation}}
\newcommand{\ba}{\begin{array}}
\newcommand{\ea}{\end{array}}

\newcommand{\bpm}{\begin{pmatrix}}
\newcommand{\epm}{\end{pmatrix}}

\newcommand{\ninn}{n \in \mathbb{N}} 

\newcommand{\bd}{\begin{definition}}
\newcommand{\ed}{\end{definition}}

\newtheorem{lemma}{Lemma}[section]
\newtheorem{theorem}{Theorem}[section]

\newtheorem{proposition}{Proposition}[section]

\newtheorem{definition}{Definition}[section] 
\thanks{$^*$ School of Mathematics, Georgia Institute of Technology, 686 Cherry Street, Atlanta, GA 30332, USA. ({\tt wong@math.gatech.edu})}
\begin{document}

\title[Asymptotic analysis of OP via the transfer matrix approach]
{Asymptotic analysis of orthogonal polynomials via the transfer matrix approach}
\author[M.-W. L. Wong]{Manwah Lilian Wong*}
\date{\today}

\maketitle

\begin{abstract}  In this paper, we present a new method via the transfer matrix approach to obtain asymptotic formulae of orthogonal polynomials with asymptotically identical coefficients of bounded variation. We make use of the hyperbolicity of the recurrence matrices and employ Kooman's Theorem to diagonalize them simultaneously. The method introduced in this paper allows one to consider products of matrices such that entries of consecutive matrices are of bounded variation. 

Finally, we apply the asymptotic formulae obtained to solve the point mass problem on the real line when the measure is essentially supported on an interval. We prove that if a point mass is added to such a measure outside its essential support, then the perturbed recurrence coefficients will also be asymptotically identical with the same limit and of bounded variation.

\end{abstract}

\section{Introduction} Let $\mu$ be a non-trivial measure on $\mathbb{R}$ such that for all $n \in \mathbb{N}$,
\be
\ds \int_{\mathbb{R}} |x|^n d\mu(x) < \infty .
\ee

We form an inner product and a norm on $L^2(\mathbb{R}, d\mu)$ as follows: for any $f,g \in L^2(\mathbb{R}, d\mu)$, define
\begin{eqnarray}
\langle f,g \rangle & = & \ds \int_{\mathbb{R}} \ol{f(x)} g(x) d\mu(x)   , \\
\|f\|^2 & = & \ds \int_{\mathbb{R}} f(x)^2 d\mu(x) .
\end{eqnarray}

By the Gram--Schmidt process, we can orthogonalize $1,x, x^2, \dots$ and obtain the family of orthogonal polynomials. We denote the $n$-th monic orthogonal polynomial as $P_n(x)$ and the $n$-th orthonormal polynomial as $p_n(x)$. Let
\be
\kappa_n = \ds \frac{1}{\|P_n\|} .
\label{kappan}
\ee Then the $n$--th orthonormal polynomial is given by
\be
p_n(x) = \ds \frac{P_n(x)}{\|P_n\|} = \kappa_n x^n + \text{ lower order terms}.
\ee

It is well-known that these orthogonal polynomials satisfy the following three-term recurrence relations:
\begin{eqnarray}
xP_n(x) & = &P_{n+1}(x) + b_{n+1} P_n(x) + a_n^2 P_{n-1}(x) \label{recurrence1}\\
x p_n(x) & = & a_{n+1}(x) p_{n+1}(x)+ b_{n+1} p_n(x) + a_n p_{n-1}(x) \label{recurrence2}
\end{eqnarray} with the properties that
\be
a_n  =  \ds \frac{\|P_n\|}{\|P_{n-1}\|} = \ds \frac{\kappa_{n-1}}{\kappa_n} > 0 \quad \text{and} \quad b_{n+1}  = \langle xp_n, p_n \rangle .
\ee

The $(a_n, b_n)_{n=1}^{\infty}$ are called the recurrence coefficients associated to the measure $d\mu$. The recurrence relation \eqref{recurrence2} is often represented by the matrix
\be
J = \bpm b_1  & a_1 & 0 & 0 & \dots \\
a_1 & b_2 & a_2 & 0 & \dots \\
0 & a_2 & b_3 & a_3 & \dots \\
\dots & \dots & \dots & \dots & \dots 
 \epm,
\label{jacobi}
\ee which is called the \textbf{Jacobi matrix}.

Apart from the Jacobi matrix, there is another representation of the recurrence relation by means of the \textbf{transfer matrix} $T_n(x)$: observe that the recurrence relation \eqref{recurrence2} can be rewritten as
\be
a_{n+1} p_{n+1}(x) = (x - b_{n+1}) p_n (x)- a_n p_{n-1} (x)\text{ for } n \geq 0 .
\label{e5}
\ee

Therefore, the general solution of \eqref{e5} above can be expressed in the following way:
\be
\bpm p_{n+1}(x) \\ a_n p_n(x) \epm = A_{n+1} (x)  \bpm p_n(x) \\ a_n p_{n-1} (x)\epm \, ,
\ee where
\be
A_j(x) = {a_j}^{-1} \bpm x-b_j & -1 \\ a_j^2 & 0 \epm \, .
\label{Ajdef}
\ee 

This motivates the definition of the \textbf{transfer matrix}
\be
T_n(x) = A_{n}(x) A_{n-1}(x) \cdots A_1(x)  \quad \text{ for } n \geq 1.
\label{Tndef}
\ee In particular, the recurrence relation could be expressed in terms of the transfer matrix applied to $(1,0)^{T}$:
\be
\bpm p_n(x) \\ a_{n-1} p_{n-1}(x) \epm  = T_n(x) \bpm 1 \\ 0 \epm .
\ee

Moreover, since $\det A_j(x) = 1$ for all $j \geq 1$,
\be
\det T_n(x) = \prod_{j=1}^{n} \det A_j(x) = 1.
\ee 

The transfer matrix $T_n(x)$ and the $A_n(x)$'s will serve as important tools when we derive the asymptotic formulae for the $p_n(x)$'s in Section \ref{asympsection}.

The reader should be reminded that besides the definitions given for $T_n(x)$ and $A_n(x)$ in \eqref{Tndef} and \eqref{Ajdef}, there are a few other commonly used definitions in the literature. 

Let $J_n$ be the truncated $n \times n$ matrix obtained from the first $n$ rows and columns of the Jacobi matrix $J$. Note that by the recurrence relation,
\be
(J_n - x ) \bpm p_0(x) \\ p_1(x) \\ \vdots \\ p_{n-2}(x) \\p_{n-1}(x) \epm = \bpm 0 \\ 0 \\ \vdots \\ 0 \\ - a_n p_n(x) \epm .
\ee

Every zero of $p_n(x)$ is an eigenvalue of $J_n$ and the orthogonal polynomials form the eigenfunction for the truncated Jacobi matrix. In other words, $p_n(x)$ is the characteristic polynomial of $J_n$. Moreover, it is known that the spectral theory of one-dimensional operators (eg., the one-dimensional discrete Schr\"odinger operator) has a lot in common with the classical theory of orthogonal polynomials on the real line. Many results were proven for the Jacobi operator and the Schr\"odinger operator in parallel (see for example \cite{birthday}).  Therefore, orthogonal polynomials have  gained the attention of both the spectral theory community and the Schr\" odinger operator community in recent years.

For a more comprehensive introduction to the theory of orthogonal polynomials on the real line, the reader may refer to \cite{chihara, birthday, nevai1, totik}.

\section{Results}

First, we consider a measure in the \textbf{Nevai class} $\mathcal{M}(a,b)$, which consists of measures with recurrence coefficients satisfying $a_n \to a, b_n \to b$. It is well-known that measures in $\mathcal{M}(a,b)$ have essential support $[b-2a, b+2a]$ (see \cite{nevai1} for a detailed discussion). In this paper, we limit ourselves to $a_n \to a \not = 0$ (Dombrowski \cite{dombrowski} showed that a Jacobi matrix with $\liminf |a_n|=0$ has empty a.c. spectrum). For the importance of the Nevai class and references to the many investigations thereof, the reader may refer to \cite{lubinsky1, mnt, nevai1, nevai2}.

Apart from the asymptotic formulae listed in Theorem \ref{asymptoticthm} below, it is interesting to note the method developed in this paper by means of applying Kooman's Theorem to the transfer matrix (see Section \ref{asympsection}). We make use of the fact that if $x_0$ is outside the essential support of the measure, then $A_n(x_0)$ is hyperbolic for all large $n$. This fact allows us to apply Kooman's theorem to simultaneously diagonalize $A_n$ and deduce asymptotic formulae for $p_n(x)$ outside the essential support. Due to the length of the proof, a sketch of the proof is provided in Section \ref{structure}.

The first result can be summarized as follows:
\begin{theorem} Let $\mu$ be a measure in $\mathcal{M}(a,b)$ with recurrence coefficients of bounded variation and $x_0 \in \mathbb{R} \backslash [b-2a, b+2a]$. Then the asymptotic formulae for $p_n(x_0)$ are as follows:
\begin{enumerate}
\item $\mu(x_0) > 0$ if and only if given any $\epsilon>0$, there exists a constant $C_\epsilon$ such that 
\be |p_n(x_0)| \leq \ds  C_\epsilon (\lambda^- + \epsilon)^n.
\ee where $|\lambda^-|<1$ is the eigenvalue of $A_\infty(x_0)$ (see \eqref{Ainftydef} for the definition of $A_\infty$).
\item $\mu(x_0) = 0$ if and only if for every $n \in \mathbb{N}$, $p_n(x_0)$ is in the form
\be
p_n(x_0)= \left( \ds \prod_{j=1}^{n} \lambda_j^+  \right) k_{n} ,
\ee where $(k_n)_n$ is a convergent sequence of bounded variation that varies according to initial conditions, and $\lambda_j^+$ is the eigenvalue of the recurrence matrix $A_j(x)$ whose norm is great than $1$ (see \eqref{Ajdef} for the definition of the recurrence matrix $A_j(x)$) . The sequence will be computed explicitly in the proof in Section \ref{asympsection}.
\end{enumerate}
\label{asymptoticthm}
\end{theorem}

Then we apply the asymptotic formulae obtained to solve the point mass problem. We add a pure point $x_0 \in \mathbb{R}$ to $\mu$ to form the measure $\tilde{\mu}$ as follows:
\be
\tilde{\mu} = \mu + \gamma \delta_{x_0} , \quad \gamma>0 .
\label{addpointreal}
\ee
 
In Theorem \ref{theorem1} below, we give formulae relating the orthogonal polynomials and the recurrence coefficients of $\mu$ and $\tilde{\mu}$. Even though these formulae are known (see \cite{nevai1}), the proofs are given below for the convenience of the reader.

\textbf{A note on notation.} We shall denote objects associated to the measure $\tilde{\mu}$ with a $\sim$ on top. For example, the $n$-th monic orthogonal polynomial with respect to the measure $\tilde{\mu}$ is denoted as $\tilde{P}_n(x)$.

\begin{theorem} Let $\mu$ and $\tilde{\mu}$ be measures defined as in \eqref{addpointreal}. Then the $n$-th monic orthogonal polynomials of $\mu$ and $\tilde{\mu}$ are related by the following formula:
\begin{eqnarray}
\tilde{P}_n(x)  &= \left( \ds \frac{\kappa_n}{\tilde{\kappa}_n}\right)^2 \left[P_n(x) -  \ds \frac{\gamma P_n (x_0) K_n(x, x_0) }{1+\gamma K_n(x_0, x_0)}\right] \, ,
\label{Pn} 
\end{eqnarray} with
\be
\left( \ds \frac{\kappa_n }{\tilde{\kappa}_n} \right)^2 =  \ds \frac{1+\gamma K_{n}(x_0, x_0)}{1 + \gamma K_{n-1}(x_0, x_0)}  ,
\ee and the reproducing kernel $K_n(x,y)$ is defined as
\be
K_n(x,y) = \ds \sum_{j=0}^{n} p_j(x) p_j(y) .
\label{Kndef}
\ee

Furthermore, the recurrence coefficients of $\mu$ and $\tilde{\mu}$ are related as follows:
\begin{enumerate}
\item \be
\tilde{a}_n = a_n \sqrt{\ds \frac{t_{n-1}}{t_n}} \, > 0
\label{antilde}
\ee where \be
t_n =  \ds \frac{1+\gamma K_{n-1}(x_0,x_0)}{1+\gamma K_{n}(x_0,x_0)}.
\label{tndef}
\ee 
\item
\be
\tilde{b}_{n+1} = b_{n+1} -  \ds \frac{\gamma P_n(x_0) p_{n-1}(x_0) \kappa_{n-1}}{1+\gamma K_{n-1}(x_0,x_0)} +   \ds \frac{\gamma P_{n+1}(x_0) p_{n}(x_0) \kappa_{n}}{1+\gamma K_{n}(x_0,x_0)}.
\label{tbn}
\ee
\end{enumerate}
\label{theorem1}
\end{theorem}

In Section \ref{prooftheorem2} we combine those asymptotic formulae with Theorem \ref{theorem1} to prove the following result:
\begin{theorem} Let $\mu$ be a non-trivial measure on $\mathbb{R}$ such that its recurrence coefficients satisfy
\begin{eqnarray}
a_n \to a \not = 0, \quad b_n \to b ; \label{condition1}\\
\ds \sum_{n=0}^\infty |a_{n+1} - a_n| + |b_{n+1} - b_n|  < \infty . \label{condition2}
\end{eqnarray}

The essential support of $\mu$ is $[b-2a,b+2a]$. If we add finitely many distinct pure points $x_1, \dots, x_k  \in \mathbb{R} \backslash [b-2a,b+2a]$ to $\mu$ as follows
\be
\mu_k = \mu + \ds \sum_{j=1}^{k} \gamma_j \delta_{x_j}, \quad \gamma_j > 0,
\label{addmanypoints}
\ee then the recurrence coefficients of $\mu_k$ satisfy \eqref{condition1} and \eqref{condition2}.
\label{theorem2}
\end{theorem}
\section{asymptotic analysis and the Point Mass Problem}\label{sectionpoint}

The point mass problem has a very long history (see the Introduction of \cite{wong1} for details) and it has its physical significance. As noted in the Introduction, results were often proven for the Jacobi operator and the Schr\"odinger operator in parallel. As a result, the point mass problem has been investigated by both the orthogonal polynomials and the mathematical physics communities.

The earliest work related to the point mass problem could be due to Wigner-von Neumann \cite{wigner}, where they constructed a potential with an embedded eigenvalue. In 1946, Borg \cite{borg} proved a well-known result concerning the Sturm--Liouville problem, which implies that if two Sturm-Liouville operators have spectra differing by a finite number of eigenvalues, then their corresponding potential functions might not be the same. Later, Gel'fand--Levitan \cite{gelfand} showed that in order to recover the potential one also needs the norming constants, which correspond to the weights of pure points in the context of orthogonal polynomials.

The point mass problem has been considered under various settings. In \cite{szwarc} Szwarc considered a measure with bounded support $S \subset [0,+\infty)$ of which the recurrence coefficients satisfy $a_n a_{n+1}^{-1} \to 1$ and $a_n/b_n \to \sqrt{A}$ for some $A \geq 0$. It was proven that if a point mass is added to the measure, then $a_n - a_n' \to 0$ and $b_n - b_n' \to 0$ as $n \to \infty$, where $a_n', b_n'$ are the  perturbed coefficients, though the specific rate of convergence was not shown.

Szwarc's paper also discusses the growth of orthogonal polynomials for measures supported on $[0,1]$ with a finite number of pure points, with $a_n \to 1$ and $b_n \to 1/2$. It was proven under such conditions,
\be
 \limsup_n |p_n(x)|^{1/n} \leq 1
\ee uniformly on $[0,1]$. Unfortunately, such results (mostly derived from \cite{ntz}), like many existing asymptotic results on $p_n(x)$, are insufficient to prove Theorem \ref{theorem1}.

In \cite{mnt}, M\'at\'e--Nevai--Totik considered a measure in the Nevai class with coefficients satisfying $a_n \to 1, b_n\to 0$ and
\be
\ds \sum_n |a_{n+1} - a_n| + |b_{n+1} - b_n| < \infty.
\ee This measure is supported on $[-1,1]$. This is a special case among those considered in Theorem \ref{theorem2}. The authors gave asymptotic formulae for $p_n(x)$ on a compact set $K \subset \mathbb{C}$ with $K \cap [-1,1]= \emptyset$ as well as on a compact set $K' \subset \supp(\mu) \backslash [-1,1]$. The authors also presented ratio asymptotic results concerning $p_n(x)/p_{n-1}(x)$. However, in both cases, we need more precise bounds on error terms to establish Theorem \ref{theorem2} (see Section \ref{asympsection}).

There is an analog of the point mass problem on the unit circle, but the problem is different in nature as the recurrence coefficients of a non-trivial probability measure on $\T$ form a one-parameter family $(\alpha_n)_{n=0}^{\infty}$ with $\alpha_n \in \mathbb{D}$ (also known as the Verblunsky coefficients in \cite{simon1, simon2}).

The problem on the unit circle reads as follows: let $\nu$ be a non-trivial probability measure on $\T=\{z \in \mathbb{C}: |x|=1\}$. One adds a pure point $\zeta \in \T$ to $\nu$ to form the new probability measure
\be
\tilde{\nu} = (1-\gamma) \nu + \gamma \delta_{\zeta} , \quad 1> \gamma>0 .
\label{dnu1def}
\ee Then one has the following classic result relating the original and the perturbed orthogonal polynomials, which is an analog of \eqref{Pn} in Theorem \ref{theorem1}.
\begin{theorem}(Geronimus \cite{geronimus1, geronimus2}) Suppose the probability measure $\tilde{\nu}$ is defined as in (\ref{dnu1def}). Then the $n$--th monic orthogonal polynomial of $\tilde{\nu}$ is given by
\be \Phi_{n}(z, \tilde{\nu}) = \Phi_n(z) - \ds \frac{\vp_n(z) K_{n-1}(z, \zeta)}{(1-\gamma) \gamma^{-1} + K_{n-1}(\zeta, \zeta)}
 \label{geronimus}
\ee where $\Phi_n(z)$ and $\vp_n(z)$ are the monic and normalized orthogonal polynomials of $\nu$ respectively; and $K_{n}(z,\zeta)  = \sum_{j=0}^{n} \ol{\vp_j(\zeta)} \vp_j(z)$.
\label{geronimustheorem}
\end{theorem}

The point mass problem on the unit circle was further investigated by Cachafeiro--Marcell\'an \cite{cm1} and Simon (see Chapter 10.13 of \cite{simon2}) using very different approaches. In \cite{wong2} Wong applied the Christoffel--Darboux formula to Simon's result and proved the point mass formula, which shows that the recurrence coefficients $(\alpha_n)_{n=0}^{\infty}$ of $d\nu$ and $d\tilde{\nu}$ are related as $\alpha_n(d\tilde{\nu}) = \alpha_n(d\nu) + \Delta_n(\zeta)$,  where
\be
\Delta_n(\zeta) = \ds \frac{(1-|\alpha_n|^2)^{1/2} \ol{\vp_{n+1}(\zeta)}\vp_n^*(\zeta)}{(1-\gamma)\gamma^{-1} + K_n(\zeta,\zeta)}; \quad K_n(\zeta,\zeta) = \ds \sum_{j=0}^{n} |\vp_j(\zeta)|^2 .
\label{deltandef}
\ee The reader may compare \eqref{deltandef} with the formulae for $\tilde{a}_n$ and $\tilde{b}_n$ in Theorem \ref{theorem1}.

In \cite{wong1, wong3, wong4} Wong studied the point mass problem on the unit circle. In particular, the class $W_p$ consisting of measures of which the Verblunsky coefficients are of $p$-generalized bounded variation was identified. It was proven that upon adding a point mass to a measure in $W_p$, we obtain a new probability measure in $W_{p+1}$. Inductively, we can add a finite number $k$ of distinct pure points to a measure in $W_p$ one after another and we will end up with a new probability measure in $W_{p+k}$. These results on the unit circle could be read in parallel to Theorem \ref{theorem2}.

\section{The Transfer Matrix Approach} \label{asympsection}

\subsection{Structure of this section} \label{structure} In Section \ref{subsection1} we consider the recurrence relation in matrix form. Since $x_0$ is not in the support of the measure, $A_j(x_0)$ is hyperbolic for all large $j$. Because of the hyperbolicity and the fact that $A_j \to A_\infty$, we can apply Kooman's Theorem to prove that there is an analytic function $U$ on a neighborhood of $A_\infty$ such that for all large $j$,  $A_j$ can be simultaneously diagonalized (see Section \ref{subsection2}).

In Section \ref{subsection3}, we consider a representation that stems out from the diagonalization.  In Proposition \ref{triprop}), we show a trichotomy based on that representation and prove Theorem \ref{asymptoticthm} for two of the three cases, which are special cases that are easier to handle.

The remaining case, being the most difficult one, will be treated in Section \ref{subsection4}. We prove several estimates for this particular case. In Section \ref{subsection5}, we combine those estimates to obtain asymptotic formulae for $p_n(x_0)$.

A diagram summarizing the results is provided at the end of the section for the convenience of the reader.

\vspace{1cm}

\textbf{From now on, let $\mu$ be a measure in $\mathcal{M}(b,a)$ with recurrence coefficients satisfying \eqref{condition1} and \eqref{condition2}.}

\subsection{The recurrence relation and $p_n(x_0)$} \label{subsection1} Note that under the conditions that $a_n \to a \not = 0$ and $b_n \to b$,
\be
\lim_{j\to\infty} A_j(x) = A_\infty(x) := a^{-1} \bpm x-b & -1 \\ a^2 & 0\epm .
\label{Ainftydef}
\ee 

$A_\infty(x)$ has eigenvalues 
\be
\lambda^{\pm} 
= \ds \frac{(x-b)\pm \sqrt{(x-b)^2 - 4a^2}}{2a} .
\label{lambdadef} \ee 
Therefore, $A_\infty(x)$ has distinct eigenvalues in $\mathbb{R}$ if and only if $x \in \mathbb{R}\backslash[b-2a, b+2a]$. In that case, one of which has absolute value strictly greater than $1$ and the other strictly less than $1$. We say that $A_\infty(x)$ is \textbf{hyperbolic}.

Now consider a fixed point $x_0 \in \mathbb{R} \backslash[b-2a, b+2a]$. Without loss of generality, we assume $x_0>b+2a$. By \eqref{lambdadef}, we have
\be
\lambda_n^+ > 1 > \lambda_n^- 
\ee for all large $n$.

From now on we will write $A_j(x_0)$ as $A_j$ and similarly for other objects that appear in the proof.

\subsection{Hyperbolicity of $A_j$ and Kooman's Theorem} \label{subsection2} Since $A_j \to A_\infty$, $A_j$ is also hyperbolic for all large $j$. This allows us to use the following result by Kooman. Adapted to suit the context of this paper, the theorem reads as follows:
\begin{theorem}[Kooman \cite{kooman1, kooman2}] Let $A$ be an $\ell \times \ell$ matrix with distinct eigenvalues. Then there exist $\epsilon > 0$ and analytic functions $U(B)$ and D(B) defined on $S_{\epsilon} = \{B: \|B - A \| < \epsilon \}$ such that \\
(1) $B = U_B D_B U_B^{-1}$, $D_B$ commutes with $A$. \\ 
(2) $U_B$ is invertible for all $B \in S_\epsilon$. \\
(3) $U_A = 1$, $D_A = A$. \\
(4) By picking a basis such that $A$ is diagonal, we can have all $D_B$ diagonal with entries being the eigenvalues of $B$.
\label{kooman}
\end{theorem}
\noindent \textit{Remarks:
\begin{enumerate}
\item The formulation of Theorem \ref{kooman} is similar to Theorem 12.1.7 of \cite{simon2}, except that in \cite{simon2} the statement was meant for quasi-unitary matrices. In fact, the proof holds as long as $A$ has distinct eigenvalues.
\item Kooman's Theorem first appeared in Theorem 1.3 of \cite{kooman1}. The first application of Kooman's Theorem to orthogonal polynomials was  made by Golinskii--Nevai \cite{golinskiinevai} for the unit circle case. They proved that if the recurrence coefficients of the measure $\mu$ on $\T$ satisfy $\alpha_n \to 0$ and if $\sum_n \|A_{n+1} - A_n\| < \infty$, then the a.c. part of the measure is positive almost everywhere on $\T$.
\end{enumerate}
}

\vspace{1cm}

Let $G$ be the matrix that diagonalizes $A_\infty$, i.e.,
\be
A_\infty = G^{-1} D G ,
\ee where
\be
D = \bpm \lambda^{+} & 0 \\ 0 & \lambda^{-} \epm
\label{Ddef}
\ee and $\lambda^{\pm}$ are the eigenvalues of $A_\infty$ (defined in \eqref{lambdadef}). Pick a basis in which $A_\infty$ is diagonal. Then by the construction of the function $D$ in Theorem \ref{kooman} above, there exists an integer
\be
N > N(\epsilon)
\label{N0def}
\ee such that $A_j$ is in some $S_\epsilon$ neighborhood of $A_\infty$ and $D_{A_j}$ is a diagonal hyperbolic matrix under this basis for all $j \geq N$. In other words, there exist diagonal matrices
\be
D_j = \bpm \lambda_j^+ & 0 \\ 0 & \lambda_j^- \epm
\label{Dndef}
\ee such that $D_{A_j(x_0)} = G D_j G^{-1} $ and the eigenvalues $\lambda_j^{\pm}$ satisfy
\be
\lambda_j^+ > 1 > \lambda_j^- ,\quad \ds \lim_{j \to \infty} \lambda_j^{\pm} = \lambda^{\pm} .
\label{lambdalimit}
\ee In fact, by \eqref{Ajdef} and a straightforward computation, we can show that $\lambda_j^{\pm}$ are roots of the characteristic polynomial
\be
z^2 - (x_0 - b_j)z + a_j^2 = 0 .
\label{quadratic}
\ee

Next, we define
\be
G_j = U_{A_j} G.
\ee Then for $j \geq N$, the matrix $A_j $ can be expressed as
\be
A_j = G_j D_j G_j^{-1} .
\label{Gjdiag}
\ee

\subsection{Representation of the recurrence relation} \label{subsection3} Let $E$ be an integer, we will choose $E$ later in the proof (see Proposition \ref{triprop}).

For $n \geq E \geq N$, the transfer matrix applied to the basis vector gives
\be
T_n \bpm 1 \\ 0 \epm =G_n D_n G_n^{-1} G_{n-1} D_{n-1} G_{n-1}^{-1} \cdots D_{{E}+1} G_{{E}+1}^{-1} G_{E} v_E ,
\label{Tneqn}
\ee where
\be
v_E = \bpm v_1^{(E)} \\ v_2^{(E)} \epm := D_{E} G_{E}^{-1} A_{{E}-1} \cdots A_1 \bpm 1 \\0 \epm .
\label{wdef}
\ee 

Following \eqref{Tneqn}, we consider the equation
\be
D_n G_{n}^{-1} G_{n-1} D_{n-1} G_{n-1}^{-1} \cdots D_{E+1} G_{E+1}^{-1} G_E v = L_n \bpm u_n v_1^{(E)} \\ w_{n} v_2^{(E)} \epm ,
\label{fndef}
\ee where $u_n$ and $w_{n}$ are defined implicitly by \eqref{fndef} above and\be
L_n = \ds \prod_{k=E+1}^{n} \lambda_k^+ \quad \text{ for } n > E .
\label{Lndef}
\ee

Also, we define
\be
\bpm v_n^{(1)} \\ v_n^{(2)} \epm = v(n) = \text{ L.H.S. of \eqref{fndef}} = L_n \bpm u_n v_1^{(E)} \\ w_{n} v_2^{(E)} \epm .
\label{wndef}
\ee

\begin{proposition} For any $E \geq N$, either $v_1^{(E)}$ or $v_2^{(E)}$  in \eqref{wdef}  is non-zero.
\label{nonzeroprop}
\end{proposition}

\begin{proof} Observe that by \eqref{wdef},
\be G_{N} v= G_{N} \bpm v_1^{(E)} \\v_2^{(E)} \epm =  A_{N} \cdots A_1 \bpm 1 \\ 0 \epm =  \bpm p_{N}(x_0) \\ a_{N-1} p_{N-1}(x_0) \epm .
\label{Gnv}
\ee 

Since $G_N$ is invertible and that $a_n > 0$ for all $n$, $v=0$ implies $p_{N}(x_0)= p_{N-1}(x_0)=0$. This contradicts the fact that the zeros of $p_n(x)$ and $p_{n-1}(x)$ strictly interlace for all $n$ and $x \in \mathbb{R}$.

\end{proof}

In the following proposition, we are going to identify the trichotomy about the pair $(v_n^{(1)}, v_n^{(2)})$ in \eqref{wndef}:

\begin{proposition} \label{triprop} Let $v_1^{(E)}$ and $v_2^{(E)}$ be defined as in \eqref{wdef} and $N$ be defined in \eqref{N0def}. Then one of the following is true:
\begin{enumerate}
\item For some $E \geq N$, both $v_1^{(E)}$ and $v_2^{(E)}$ are non-zero.
\item For all $E \geq N$, $v_2^{(E)} \equiv 0$. In that case,
\be
p_j(x_0) = \left( \ds \prod_{k=N}^{j-1} \lambda_k^+ \right) p_N(x_0) .
\label{expogrowth}
\ee
\item For all $E \geq N$, $v_1^{(E)} \equiv 0$. In that case,
\be
p_j(x_0) = \left( \ds \prod_{k=N}^{j-1} \lambda_k^- \right) p_N(x_0) .
\label{expodecay}
\ee
\end{enumerate}
Remark: Case (1) will be treated in Section \ref{subsection4} below.
\end{proposition}

\begin{proof} Suppose for all $E \geq N$, $v_2^{(E)} \equiv 0$. By Proposition \ref{nonzeroprop}, $v_1^{(E)} \not = 0$ for all $E$. That is equivalent to
\be
G_E \bpm v_1^{(E)} \\ 0 \epm = \bpm p_E(x_0) \\ a_{E-1} p_{E-1}(x_0) \epm , \quad\forall E \geq N.
\ee

Recall from \eqref{Gjdiag} that $A_E = G_E D_E G_E^{-1}$. In other words, $G_E$ is the change of basis matrix that maps the vector $(1,0)$ to the eigenvector of $A_E$ with eigenvalue $\lambda_E^+$. Therefore,
\be
\bpm p_{E+1}(x_0) \\ a_E p_E(x_0) \epm = A_E \bpm p_E(x_0) \\ a_{E-1} p_{E-1}(x_0) \epm = \lambda_E^+ \bpm p_E(x_0) \\ a_{E-1} p_{E-1}(x_0) \epm .
\label{eee}
\ee 

By an inductive argument, we obtain \eqref{expogrowth}.

The proof for (3) is identical in nature except that we have
\be
\bpm p_{E+1}(x_0) \\ a_E p_E(x_0) \epm  = \lambda_E^- \bpm p_E(x_0) \\ a_{E-1} p_{E-1}(x_0) \epm .
\label{eee2}
\ee instead of \eqref{eee} because $G_E$ maps $(0,1)$ to the eigenvector of $A_E$ with eigenvalue $\lambda_E^-$. This proves \eqref{expodecay}.

\end{proof}

Now that we have asymptotic formulae for $p_n(x_0)$ for Cases (2) and (3), we are going to consider Case (1). In fact, if there exists $E$ such that $v_1^{(E)}$ and $v_2^{(E)}$ are non-zero, the asymptotic formula for $p_n(x_0)$ will be quite similar to \eqref{expogrowth} or \eqref{expodecay} except for the presence of error terms, which will be analyzed in the next section.

\vspace{1cm}

\subsection{Several estimates} \label{subsection4} By Proposition \ref{triprop}, there are three cases. In this section we focus on Case (3): i.e., there exists an $E \geq N$ such that $v_1^{(E)}$ and $v_2^{(E)}$ are both non-zero.

Consider a fixed $E$. For convenience, we shall write 
\be
v_1 = v_1^{(E)}, v_2 = v_2^{(E)} \text{ and } v= v_E.
\label{v1v2}
\ee

\begin{proposition} \label{prop1} There is a constant $C$ such that
\be
\|v(n+1) - D_{n+1} v(n) \| \leq C \|A_{n+1} - A_n \| |L_n| \left( |u_n| + |w_{n}|\right).
\label{bound3}
\ee 
\end{proposition}

\begin{proof} Note that
\be
v(n+1) - D_{n+1} v(n) = D_{n+1} \left( G_{n+1}^{-1} G_{n} - 1 \right) v(n).
\label{wnbound}
\ee

The goal is to give bounds for each of the components on the right hand side of (\ref{wnbound}). Since $U$ is analytic on $S_\epsilon$, on some compact subset of $S_\epsilon$ there exist constants $\eta_1, \eta_2>0$ such that
\be
\|G_{n+1}- G_{n}\|  \leq \|G\| \|U_{A_{n+1}} - U_{A_{n}}\|  \leq \eta_1 \| A_{n+1} - A_{n}\| \label{bound1} \ee and
\be
\| G_{n+1}^{-1} \| \leq  \| G^{-1}\| \| U_{A_{n+1}}^{-1}\| \leq  \eta_2 . \label{bound2}
\ee

Thus, for $\eta = \eta_1 \eta_2$,
\be
\|G_{n+1}^{-1} G_{n} - 1\| = \| G_{n+1}^{-1} \left(  G_{n}-G_{n+1} \right) \| \leq \eta \| A_{n+1} - A_n\| .
\label{bound1a}
\ee

Furthermore,
\be
\ds \sup_{n \geq N} \| D_{n}\| = \sup_{n \geq N} |\lambda_{n}^+| < 2 |\lambda^+|
\label{bound1b}
\ee and 
\be
\|v(n)\|  = \left\| \bpm u_n L_n v_1 \\ w_{n} L_n v_2 \epm \right\| < C_1 |L_n| \left( |u_n| + |w_{n}|\right) ,
\label{bound1c}
\ee where $C_1=max\{|v_1|, |v_2|\} > 0$. By applying \eqref{bound1a}, \eqref{bound1b} and \eqref{bound1c} to (\ref{wnbound}), we finish the proof of Proposition \ref{prop1}.
\end{proof}

\vspace{1cm}

\begin{proposition} \label{wnto0} Let $u_n$ and $w_{n}$ be defined as in \eqref{fndef} and $v_1, v_2 \not = 0$. Then the following inequalities hold:
\begin{enumerate}
\item There is a constant $C_3$ such that
\be \left| u_{n+1} - u_n \right|  \leq C_3  \| A_{n+1} - A_n\|  \left( |u_n| + |w_{n}|\right) .
\label{bound6}
\ee
\item There is a constant $C_4$ such that
\be 
\left| w_{n+1} - \ds \frac{\lambda_{n+1}^-}{\lambda_{n+1}^+} w_{n}\right|  \leq C_4 \| A_{n+1} - A_n\|  \left( |u_n| + |w_{n}|\right) .
\label{bound7}
\ee
\end{enumerate}
\end{proposition}

\begin{proof}

Recall that $L_{n+1} =\lambda_{n+1}^+ L_n$ and $v_n^{(1)}=L_n u_n v_1$. Since $v_1 \not = 0$,
\begin{multline}
 \left| u_{n+1} - u_n \right|  = \ds   \left| \ds \frac{v_{n+1}^1 -\lambda_{n+1}^+ v_n^{(1)}}{ v_1 L_{n+1}} \right|
\leq \ds \frac{\|w(n+1) - D_{n+1} w(n)\|}{| v_1 L_{n+1}|}.
\label{v1}
\end{multline}

Similarly, since $v_n^{(2)} = L_n w_n v_2$ and $v_2 \not = 0$,
\be
\left|w_{n+1} -  \ds \frac{\lambda_{n+1}^-}{\lambda_{n+1}^+} w_n \right| = \left| \ds \frac{v_{n+1}^{(2)} - \lambda_{n+1}^- v_n^{(2)} }{v_2 L_{n+1}} \right| \leq \ds \frac{\|w(n+1) - D_{n+1} w(n)\|}{|v_2 L_{n+1} |}.
\label{v2}
\ee

Apply Proposition \ref{prop1} to the equations above to obtain \eqref{bound6} and \eqref{bound7}.

\end{proof}

The following lemma concerning $u_n$ and $w_{n}$ is central for this paper. As we shall see in the proof of Theorem \ref{theorem2}, it implies the dichotomy between exponential decay and exponential growth of $p_n(x_0)$.

\begin{lemma} Let $u_n$ and $w_{n}$ be defined as in \eqref{fndef} and $v_1, v_2 \not = 0$. Then one of the following is true:
\begin{enumerate}
\item There exists a constant $C$ such that $|u_n| \leq C |w_{n}|$. Moreover, given any $\epsilon>0$, there exist an integer $N_{\epsilon}$ and a constant $C_\epsilon$ such that
\be
|w_{n}| \leq C_\epsilon \left( \ds \left|\frac{\lambda^-}{\lambda^+} \right|+ \epsilon \right)^{n} , \quad \forall n \geq N_\epsilon.
\ee 
\item $|w_{n}/ u_n| \to 0$. Furthermore, ${u_\infty} = \lim_{n \to \infty} u_n$ exists and it is non-zero.
\end{enumerate}
\label{lemmao1}
\end{lemma}

\begin{proof} There are two possible situations concerning $u_n$ and $w_{n}$:
\begin{enumerate}
\item There exist a fixed integer $K_0$ and a constant $C$ such that $|u_n| \leq C |w_{n}|$ for all $n \geq K_0$.
\item For any integer $K$ and any constant $H$, there exists an integer $n_{K,H} \geq K$ such that 
$|u_{n_{K,H}}| > H |w_{n_{K,H}}|$.
\end{enumerate}

\noindent \textbf{Case (1)}: By (\ref{bound7}), for $n \geq \max\{N, K_0 \}$, there is a positive constant $C_7$ such that
\be
|w_{n+1}| \leq \left( \left| \ds \frac{\lambda_{n}^-}{\lambda_{n}^+}\right| + C_7 \|A_{n+1} - A_n \| \right) |w_{n}|.
\ee

Recall that $\|A_{n+1} - A_n\| \to 0$ and $\lambda_{n}^{\pm} \to \lambda^{\pm}$. Thus, given any $\epsilon>0$, there exist an integer $N_{\epsilon}$ and a constant $C_\epsilon$ such that
\be
|w_{n}| \leq C_\epsilon \left( \ds \left|\frac{\lambda^-}{\lambda^+} \right|+ \epsilon \right)^{n} \quad \forall n \geq N_\epsilon .
\ee In other words, $w_{n}$ decays exponentially fast. Hence, $u_n$ also decays exponentially fast to zero given the inequality $|u_n| \leq C|w_n|$. This corresponds to (2a) of Lemma \ref{lemmao1}.

\vspace{1cm}

\noindent \textbf{Case (2)}: Let $r_n = w_{n}/u_n$. The statement $r_n \to 0$ is by definition equivalent to proving that given any $\epsilon >0$ there exists an integer $J_{\epsilon}$ such that $|r_j| < \epsilon$ for all $j \geq J_{\epsilon}$.

First, we show that both $u_n$ and $u_{n+1}$ are non-zero, as \eqref{eqn16} below will involve $u_n$ and $u_{n+1}$ in the denominator.

By the assumption we can choose any $H$, we choose one such that $1/H < \epsilon$. Consider any fixed pair $(K, H)$ (the choice of $K$ will be made later in the proof). The existence of an integer $n = n_{K,H} > K$ such that $|r_n| < 1/H=\epsilon$ is guaranteed, which implies $u_n \not = 0$. Furthermore, by the triangle inequality and (\ref{bound6}),
\be \begin{array} {lll}
\left| \ds \frac{u_{n+1}}{u_n} \right| & \geq & 1- \left|\ds \frac{u_{n+1} - u_n}{u_n} \right| \\
\\
& \geq & 1 - C_3 \|A_{n+1} - A_n \| (1+ |r_n| ) > 0 
\end{array}
\label{eqn12}
\ee which implies that $u_{n+1}$ is also non-zero.

Next, observe that
\be \begin{array}{lll}
&\left| r_{n+1} - \ds \frac{\lambda_{n}^-}{\lambda_{n}^+} r_n\right|\\
\leq & \left| \ds \frac{w_{n+1}}{u_{n+1}} - \frac{\lambda_{n}^-}{\lambda_{n}^+} \frac{w_{n}}{u_{n+1}} \right| + \ds \left| \ds \frac{\lambda_{n}^-}{\lambda_{n}^+}\right| \left| \ds \frac{w_{n}}{u_{n+1}} - \ds \frac{w_{n}}{u_n}\right| \\
\\
=&  \left| \ds \frac{w_{n+1} - (\lambda_{n}^-/\lambda_{n}^+) w_{n}}{u_{n+1}}\right| + \left| \ds \frac{\lambda_{n}^-}{\lambda_{n}^+} r_n \right| \left| \ds \frac{u_n-u_{n+1}}{u_{n+1}}  \right| .
\end{array}
\label{eqn16}
\ee

By (\ref{bound6}) and ({\ref{bound7}}), there is a positive constant $C_8$ such that
\be \begin{array}{lll}
& \left| r_{n+1} - \ds \frac{\lambda_{n}^-}{\lambda_{n}^+} r_n\right| \\
\\
\leq & \ds \frac{ 1+|r_n||\lambda_{n}^-/\lambda_{n}^+|}{|u_{n+1}|} C_8  \| A_{n+1} - A_n\| (|u_n|+|w_{n}| )\\
=&  C_8 (1+|r_n||\lambda_{n}^-/\lambda_{n}^+|) \|A_{n+1} - A_n \|  \ds \frac{  |u_n|}{|u_{n+1}|} (1+ |r_n| ) .
\label{eqn11}
\end{array}
\ee

By inverting (\ref{eqn12}),
\be
\left| \ds \frac{u_n}{u_{n+1}} \right| \leq \ds \frac{1}{1 - C_3\|A_{n+1} - A_n \| (1+ |r_n| )} \, .
\ee

Then we plug it into (\ref{eqn11}) to get
\be
\left| r_{n+1}\right|  \leq \left| \ds \frac{\lambda_{n}^-}{\lambda_{n}^+} r_n\right| +
 \ds \frac{C_8(1+|r_n||\lambda_{n}^-/\lambda_{n}^+|)  (1+ |r_n| )}{1 - C_3\|A_{n+1} - A_n \| (1+ |r_n| )}  \|A_{n+1} - A_n \| .
\label{eqn13}
\ee


Since $\|A_{n+1} - A_n \| \to 0$, the second term on the right hand side of (\ref{eqn13}) can be arbitrarily small if $n$ is sufficiently large. Hence, for any sufficiently large $K$, there exists $n > K$ such that $|r_{n+1}|<|r_n|<\epsilon$.

Inductively, we can aply the same argument to $r_{n+1}$ to prove that $|r_{n+2}| < \epsilon$. As a result, $|r_{j}| < \epsilon$ for all large $j$. This proves $|w_{n}/u_n| \to 0$.\\

Next, we are going to show that $\lim_{n \to \infty} u_n$ exists. Divide (\ref{bound6}) by $|u_n|$. Since $|r_n| \to 0$,
\be
\left| \ds \frac{u_{n+1}}{u_n} - 1 \right| \leq C  \|A_{n+1} - A_n\| \left( 1 + |r_n| \right) \to 0 .
\label{eqn18a}
\ee

Moreover, since $\log$ is analytic near $1$, in an $\epsilon$-neighborhood of $1$ there is a constant $E$ such that
\be |\log z| = |\log \zeta - \log 1| \leq E |z-1| .
\ee

By (\ref{eqn18a}),
\be
\left| \log \left(\ds \frac{u_{n+1}}{u_n} \right) \right| \leq C \|A_{n+1} - A_n\| .
\label{eqn19}
\ee

Therefore, the series $\sum_{j=N}^{\infty} \log \left( u_{j+1}/u_{j} \right)$ is absolutely convergent. Furthermore, as we have seen in \eqref{eqn12}, $u_{j} \not = 0$ for all large $j$. Thus, $\log u_{j}$ is well-defined and the following limit 
\be
u_\infty : = \lim_{n\to \infty} \log u_{n}  = \lim_{n \to \infty}\ds \sum_{j=p}^{n-1} \left( \log  u_{j+1} - \log u_{j} \right) + \log u_p
\ee exists and is finite.

This corresponds to the second part of (2b) and concludes the proof of Lemma {\ref{lemmao1}}.

\end{proof}

\subsection{Implications of Lemma \ref{lemmao1}} \label{subsection5}

By Lemma \ref{lemmao1}, there are two possible situations: \\

\textbf{Case (1).} Suppose (2a) of Lemma \ref{lemmao1} is true. We will prove that this corresponds to the case $\mu(x_0)>0$.

First, observe that
\be
\bpm p_{n}(x_0) \\ a_{n-1} p_{n-1}(x_0) \epm = T_n(x_0) \bpm 1 \\ 0 \epm = G_n L_n \bpm u_n v_1 \\ w_{n} v_2 \epm .
\ee

Hence, given any $\epsilon >0$, there exists a constant $K_\epsilon>0$ such that
\begin{multline}
\left\| T_n(x_0) \bpm 1 \\ 0 \epm \right\| \leq \|G_n\| \ds \prod_{k=N+1}^n |\lambda^+_k| \left\| \bpm u_n v_1 \\ w_{n} v_2 \epm \right\| \\
\leq K_\epsilon (|\lambda^+|+\epsilon)^n \left( \ds \left|\frac{\lambda^-}{\lambda^+} \right|+ \epsilon \right)^{n}.
\label{e200}
\end{multline}

Since $|\lambda^-|<1$ and $ |\lambda^-/\lambda^+| <1$, $p_{n}(x_0)$ goes to zero exponentially fast by \eqref{e200}.

It is well-known (see for example \cite{dks}) that
\be
\mu(x_0) = \left( \ds \lim_{N \to \infty} \sum_{n=0}^{N} p_n(x_0)^2 \right)^{-1} .
\label{pointsum}
\ee

Hence, $\mu(x_0)>0$ and we are just varying the weight of $x_0$. It is easy to prove that by \eqref{tndef} and \eqref{tbn}, both $a_n - \tilde{a}_{n}$ and $\tilde{b}_{n+1} -b_{n+1}$ go to $0$ exponentially fast.

This is in agreement with Simon's result (Corollary 24.4 of \cite{simon3}) that varying the weight of a pure point will result in exponentially small perturbation of the recurrence coefficients.

\vspace{1cm}

\textbf{Case (2).} Suppose (2b) of Lemma \ref{lemmao1} is true. Let
\be
G_n = \bpm g_{1,n} & g_{1,n}' \\ g_{2,n}  & g_{2,n}' \epm \to G = \bpm g_1  & g_1' \\ g_2 & g_2' \epm .
\label{gdef2}
\ee

Since $p_{n}(x_0)$ is the first component of the vector $G_n L_n (u_n v_1, w_{n} v_2)^T$, we have
\be \begin{array}{ll}
p_{n}(x_0) & = L_n (g_{1,n} u_n v_1 + g'_{1,n} w_{n} v_2 ) \\
& = L_n u_n \left(g_{1,n} + g'_{1,n} r_n v_2 \right) \\
& = L_n ({u_\infty} g_1 v_1 + o(1)) .
\label{pnasymp}
\end{array}
\ee The last equality holds because $r_n = u_n/w_{n} \to 0$ by (2b) of Lemma \ref{lemmao1}.

Here is a summary of the results in this section:
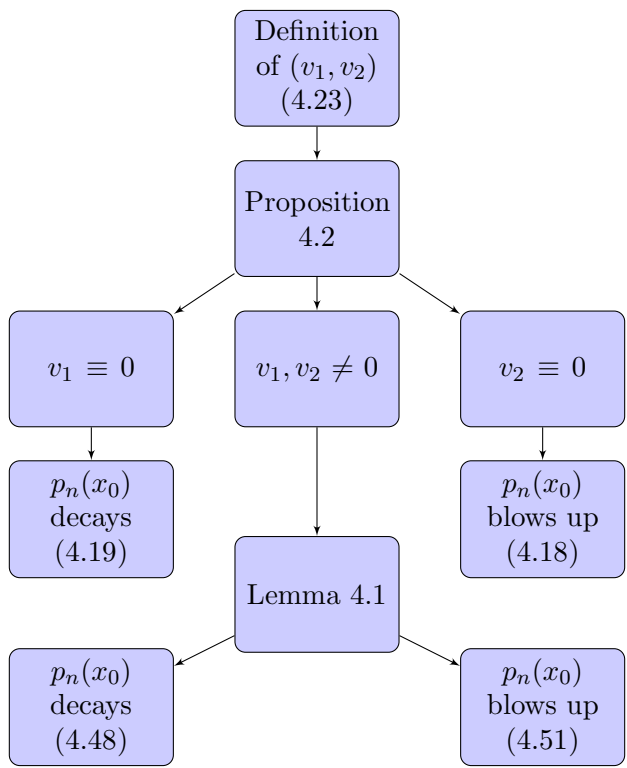
\begin{center}
\tikzstyle{decision} = [diamond, draw, fill=blue!20, 
    text width=4.5em, text badly centered, node distance=3cm, inner sep=0pt]
\tikzstyle{block} = [rectangle, draw, fill=blue!20, 
    text width=5em, text centered, rounded corners, minimum height=4em]
\tikzstyle{line} = [draw, -latex']
\tikzstyle{cloud} = [draw, ellipse,fill=red!20, node distance=3cm,
    minimum height=2em]
   \begin{tikzpicture}[node distance = 2cm, auto]
    \node [block] (def) {Definition of $(v_1, v_2)$ \\ \eqref{v1v2}};
    \node [block, below of=def, node distance=2cm] (triprop) {Proposition \ref{triprop}};
    \node [block, below of=triprop, node distance=2cm] (nonzero) {$v_1, v_2 \not = 0$};
    \node [block, left of=nonzero, node distance=3cm] (v1zero) {$v_1\equiv 0$};
    \node [block, right of=nonzero, node distance=3cm] (v2zero) {$v_2 \equiv 0$};
      \node [block, below of =v1zero, node distance = 2cm](decay1){$p_n(x_0)$ decays \\ \eqref{expodecay}};
    \node [block, below of =v2zero, node distance =2cm](blows1){$p_n(x_0)$ blows up \\ \eqref{expogrowth} };   
    \node [block, below of = nonzero, node distance = 3cm](lemma) {Lemma \ref{lemmao1}};
    \node [block, below of =v1zero, node distance = 4.5cm](decay){$p_n(x_0)$ decays \\ \eqref{e200}};
  \node [block, below of = v2zero, node distance =4.5cm](blows){$p_n(x_0)$ blows up \\ \eqref{pnasymp}};

\path [line] (def) -- (triprop);
\path [line] (triprop) --(v1zero);
\path [line] (triprop) -- (nonzero);
\path [line] (triprop) -- (v2zero);
\path [line] (nonzero) -- (lemma);
\path [line] (lemma) -- (decay);
\path [line] (lemma) -- (blows);
\path [line] (v1zero) -- (decay1);
\path [line] (v2zero) -- (blows1);
\end{tikzpicture}

\end{center}

\section{Proof of Theorem \ref{asymptoticthm}} \label{proofasymptoticthm}By the discussion in Section \ref{asympsection}, $p_n(x_0)^2$ is either exponentially increasing or exponentially decaying towards zero. Moreover, by \eqref{pointsum} above,
\be
\mu(x_0) > 0 \Longleftrightarrow \ds \sum_{n=1}^{\infty} p_n(x_0)^2 < \infty .
\ee 

Therefore, $\mu(x_0) > 0$ if and only if $p_n(x_0)^2$ is exponentially decaying towards zero, which corresponds to \eqref{expodecay} and \eqref{e200}; $\mu(x_0) = 0$ if and only if $p_n(x_0)^2$ is exponentially increasing, which corresponds to \eqref{expogrowth} and \eqref{pnasymp}.

\section{Proof of Theorem \ref{theorem1}}\label{section1}

\begin{proof} Let $K_n(x,y)$ be the reproducing kernel of the measure $\mu$, which is given by
\be
K_n(x,y) = \ds \sum_{j=0}^{n} p_j (x) p_j(y) .
\ee

Since $\tilde{P}_n(x)$ is a polynomial of degree $n$,
\be
\begin{array}{ll}
\tilde{P}_n (x) & = \ds \int \tilde{P}_n (y) K_n(x,y) d\mu(y) \\
& = \ds \int \tilde{P}_n (y) K_n(x,y) d\tilde{\mu}(y) - \gamma \tilde{P}_n (x_0) K_n(x, x_0) \, .
\end{array}
\label{e1}
\ee

Moreover, $\tilde{P}_n(y)$ is orthogonal to all polynomials with degree $\leq n-1$ with respect to the inner product $\langle \, , \, \rangle_{d\tilde{u}}$. Therefore,
\be \begin{array}{ll}
\ds \int \tilde{P}_n(y) K_n(x,y) d\tilde{\mu}(y) & = \ds \int \tilde{P}_n (y) p_n (y) p_n(x) d\tilde{\mu}(y) \\ & = p_n(x) \ds \int \tilde{P}_n(y) p_n(y) d\tilde{\mu}(y) 
\end{array}
\label{e1a}
\ee and that
\be \left\langle \tilde{P}_n(y), p_n(y) \right\rangle_{d\tilde{\mu}} = \left\langle \ds  \tilde{P}_n(y), \kappa_n y^n \right\rangle_{d\tilde{\mu}} = \frac{\kappa_n}{(\tilde{\kappa}_n)^2} \,.
\ee
Therefore, by \eqref{e1} and \eqref{e1a},
\be
\tilde{P}_n(x) = \left( \ds \frac{ \kappa_n}{\tilde{\kappa}_n} \right)^2 P_n(x) -  \gamma \tilde{P}_n (x_0) K_n(x, x_0) .
\label{e2}
\ee

Now plug in $x=x_0$ into \eqref{e2}. Upon rearranging, we get
\be
\tilde{P}_n(x_0) = \left( \ds \frac{\kappa_n}{\tilde{\kappa}_n} \right)^2 \ds \frac{P_n(x_0)}{1+\gamma K_n(x_0, x_0)} .
\label{e3}
\ee 

Putting \eqref{e3} into \eqref{e2}, we arrive at \eqref{Pn}. In particular, note that both $P_n(x)$ and $\tilde{P}_n(x)$ are monic polynomials. Hence, by comparing the coefficients of $x^n$ on each side of \eqref{Pn}, we get 
\begin{eqnarray}
1 & = & \left( \ds \frac{\kappa_n}{\tilde{\kappa}_n} \right)^2 \left( 1 - \ds \frac{\gamma p_n(x_0)^2}{1 + \gamma K_n(x_0, x_0)} \right) \label{kappafrac1}\\
& = & \left( \ds \frac{ \kappa_n}{\tilde{\kappa}_n} \right)^2 \ds \frac{1+\gamma K_{n-1}(x_0, x_0)}{1 + \gamma K_n(x_0, x_0)}  .
\label{kappafrac2}
\end{eqnarray}

Recall that $a_n = \kappa_{n-1}/\kappa_{n}$. This proves \eqref{antilde}.

Next, we prove the formula for $\tilde{b}_{n+1}$. Let $m_n$ be the coefficient of $x^{n-1}$ in $P_n(x)$, i.e.,
\be
P_n(x) = x^n + m_{n} x^{n-1} + \text{ lower order terms }.
\ee By the recurrence relation \eqref{recurrence1}, $b_{n+1}$ is given by the coefficient of $x^n$ in $xP_n(x)-P_{n+1}(x)$, which can also be expressed as
\be
b_{n+1} = m_n - m_{n+1}.
\ee

To prove formula \eqref{tbn} for $\tilde{b}_{n+1}$, we will compute $\tilde{m}_n$. By \eqref{Pn},
\be 
\tilde{m}_n = \left( \ds \frac{\kappa_n}{\tilde{\kappa}_n} \right)^2 \left(m_n - \ds \frac{\gamma P_n(x_0)}{1+ \gamma K_n(x_0,x_0)} \left[ p_n(x_0) \kappa_n m_n + p_{n-1}(x_0) \kappa_{n-1} \right] \right).
\label{e15}
\ee

The coefficients of $m_n$ in \eqref{e15} are given by
\be 
\left( \ds \frac{\kappa_n}{\tilde{\kappa}_n} \right)^2 \left(1 - \ds \frac{\gamma p_n(x_0)^2}{1+\gamma K_n(x_0,x_0)}\right) ,
\ee which is equal to $1$ by \eqref{kappafrac1}. Therefore,
\be \begin{array}{lll}
\tilde{m}_n & =m_n - \left( \ds \frac{\kappa_n}{\tilde{\kappa}_n} \right)^2 \ds \frac{\gamma P_n(x_0) p_{n-1}(x_0) \kappa_{n-1}}{1+\gamma K_n(x_0,x_0)} \\
\\
& = m_n -  \ds \frac{\gamma P_n(x_0) p_{n-1}(x_0) \kappa_{n-1}}{1+\gamma K_{n-1}(x_0,x_0)} .
\end{array}
\ee The last equality follows from the expression of $(\kappa_n/\tilde{\kappa}_n)^2$ in \eqref{kappafrac2}.

This concludes the proof of Theorem \ref{theorem1}. \end{proof}

\section{Proof of Theorem \ref{theorem1}} \label{prooftheorem2}

We are going to separate the proof into two different cases. \\

\noindent \textbf{Case (1): $x_0$ is a pure point of $\mu$.} This is the easier case. Since
\be
\ds \lim_{n \to \infty} K_n(x_0,x_0) = \mu(x_0)^{-1},
\ee it is clear by \eqref{tndef} that $t_n \to 1 $ and $\tilde{a}_n \to a$ as $n \to \infty$. Furthermore, recall that $P_{n+1}(x_0) = p_{n+1}(x_0)/ \kappa_{n+1}$ and $a_{n+1} = \kappa_n/ \kappa_{n+1}$. Hence,
\be
 h_n:= \ds \frac{\gamma P_{n+1}(x_0) p_n(x_0) \kappa_n}{1+\gamma K_n(x_0,x_0)} =\ds \frac{\gamma a_{n+1} p_{n+1}(x_0) p_n(x_0)}{1+\gamma K_n(x_0,x_0)} \to 0
\label{hndef}
\ee because $p_n(x_0) \to 0$ exponentially fast. Therefore, 
\be
\tilde{b}_{n+1} = b_{n+1} - h_{n-1} + h_n \to b \text{ as } n \to \infty .
\label{bnlimit}
\ee It is trivial to prove \eqref{condition2} for this particular case so we shall omit the proof.

\vspace{1cm}

\noindent \textbf{Case (2): $x_0$ is not a pure point of $\mu$.} To prove that $\tilde{a}_n \to a$, we are going to prove that $\lim_{n\to \infty} t_n$ exists by employing the following theorem:
\begin{theorem}[Ces\`aro--Stolz \cite{cesaro}] Let $(\Gamma_n)_{\ninn}, (K_n)_{\ninn}$ be two sequences of numbers such that $K_n$ is strictly increasing and tends to infinity. If the following limit 
\be
\ds \lim_{n \to \infty} \frac{\Gamma_n - \Gamma_{n-1}}{K_n - K_{n-1}}
\ee exists then it is equal to $\lim_{n \to \infty} \Gamma_n/K_n$.
\label{cstheorem}
\end{theorem}

Let $\Gamma_n = 1 + \gamma K_{n-1}(x_0,x_0)$ and $K_n = 1 + \gamma K_n(x_0, x_0)$. Since $\lambda_n^+ > 1$, $L_n \to \infty$. Hence, $p_n(x_0) \to \infty$, which also implies that $K_n$ is an increasing sequence to infinity. This allows us to apply the Ces\`aro--Stolz Theorem to prove that
\be
\ds \lim_{n \to \infty} t_n
=\lim_{n \to \infty}   \ds \frac{\Gamma_n - \Gamma_{n-1}}{K_n - K_{n-1}} 
= \ds \lim_{n \to \infty} \ds \frac{p_{n-1}(x_0)^2}{p_n(x_0)^2} = \ds \frac{1}{{\lambda^+}^2} .
\label{tnlimit}
\ee The last equality of \eqref{tnlimit} follows from the asymptotic formula \eqref{pnasymp} for $p_n(x_0)$. As a result, $\lim_{n \to \infty} \tilde{a}_n = \lim_{n \to \infty} a_n = a$.

\vspace{1cm}

Next, we are going to prove that $\lim_{n \to \infty} \tilde{b}_n = b$. 
We shall employ the Ces\`aro--Stolz Theorem again. Let
\be
\Gamma_n = \gamma p_{n+1}(x_0) p_n(x_0) \text{ and } K_n = 1+ \gamma K_n(x_0,x_0) .
\ee

By the asymptotic formula \eqref{pnasymp} for $p_n(x_0)$,
\be
\ds \frac{\Gamma_n - \Gamma_{n-1}}{K_n - K_{n-1}} = \ds \frac{ \gamma p_n(x_0)\left(p_{n+1}(x_0) - p_{n-1}(x_0) \right)}{\gamma p_n(x_0)^2} \to  \lambda^+ - \ds \frac{1}{\lambda^+} .
\ee 

Therefore,
\be
\ds \lim_{n \to \infty} h_n= \ds \lim_{n \to \infty} \ds \frac{\gamma a_{n+1} p_{n+1}(x_0) p_n(x_0)}{1+\gamma K_n(x_0,x_0)} = a\left( \lambda^+ - \ds \frac{1}{\lambda^+} \right) 
\ee and by \eqref{bnlimit}, $\tilde{b}_n \to b$.

\vspace{1cm}

We are going to prove that the sequence $(\tilde{a}_n)_n$ is of bounded variation.

Since $\lim_{n \to \infty} \tilde{a}_n = a \not = 0$, for all large $n$ there exists a constant $C$ such that
\be
|\tilde{a}_{n+1}^2 - \tilde{a}_n^2| = |\tilde{a}_{n+1} - \tilde{a}_n| |\tilde{a}_{n+1} +\tilde{a}_n| \geq C |\tilde{a}_{n+1} - \tilde{a}_n| .
\ee 

Therefore, it is enough to prove that $(\tilde{a}_n^2)_n$ is of bounded variation. Furthermore, by the formula for $\tilde{a}_n$ \eqref{antilde} and the fact that $\lim_{n \to \infty} t_n = 1/{\lambda^+}^2$, it suffices to show that the sequence $(t_n)_{n}$ is of bounded variation.

To do that, we are going to show that the sequence $(\tau_n)_n$ is of bounded variation, where
\be
\tau_n = \ds \frac{p_n(x_0)^2}{\gamma^{-1} + K_n(x_0, x_0)} \, .
\label{taundef}
\ee

There are two cases to consider, $p_n(x_0)$ being in the form  \eqref{expogrowth} or $p_n(x_0)$ being in the form \eqref{pnasymp}.

If $p_n(x_0)$ is in the form \eqref{expogrowth},
\be
\tau_n = \ds \frac{ L_n^2}{\gamma^{-1}+ K_n(x_0, x_0)} \ds \frac{p_N(x_0)^2}{{\lambda_n^+}^{2}}
\label{taun1}
\ee

If $p_n(x_0)$ is in the form \eqref{pnasymp},
\be
\tau_n =\underbrace{\ds \frac{L_n^2}{\gamma^{-1} + K_n(x_0,x_0)} }_{\textrm{(I)}} \underbrace{u_n^2 \left(g_{1,n} v_1 + g'_{1,n} r_n v_2 \right)^2 }_{\textrm{(II)}} .
\label{taun2}
\ee

Observe that $p_N(x_0)^2/(\lambda_N^+)^2$ is of bounded variation, because $p_N(x_0)^2$ is a constant and $(1/(\lambda_n^+)^2)_n$ is of bounded variation. Moreover, $|g_{1,n+1}- g_{1,n}|$, $|g_{1,n+1}' - g_{2,n}'|$ and $|r_{n+1} - r_n|$ are of the order $O(\|G_{n+1} - G_n\|)$, which implies they are also of bounded variation. Therefore, given \eqref{taun1} and \eqref{taun2}, it remains to show that $L_n^2/(\gamma^{-1} + K_n(x_0, x_0))$ is of bounded variation.

We will make use of the simple fact: given the equality
\be
\ds \frac{1}{y_{n}} - \frac{1}{y_{n+1}} = \ds \frac{y_{n+1} - y_n}{y_{n+1} y_n} ,
\label{e10}
\ee if $\lim_{n \to \infty}y_n = y \not = 0$, $(y_n)_{n}$ is of bounded variation if and only if $(1/y_n)_{n}$ is of bounded variation.

Hence, we will prove that $(\gamma^{-1} + K_n(x_0, x_0))/L_n^2$ is of bounded variation and its limit exist when $n$ goes to infinity. To prove the latter, observe that by the Ces\`aro--Stolz Theorem,
\be
\ds \lim_{n \to \infty} \ds \frac{L_{n+1}^2 - L_n^2}{K_{n+1}(x_0,x_0) - K_n(x_0, x_0)} = \ds \frac{L_n^2 \left((\lambda_{n+1}^+)^2 - 1 \right)}{p_{n+1}(x_0)^2} = \ell \left((\lambda^+)^2 - 1 \right) \not = 0 ,
\ee where $\ell$ is a non-zero constant whose value depends on whether $p_n(x_0)$ is in the form \eqref{expogrowth} or \eqref{pnasymp}.

For the convenience of computation we will define a few more objects below. First, we let
\be \Lambda_n = \begin{cases} \lambda_{n}^+ & \mbox{ if } n \geq N+1 \\
1 & \mbox{ if } 0 \leq n \leq N
\end{cases} .
\ee Then by (\ref{Lndef}), $L_n = \prod_{j=1}^{n} \Lambda_j$. Moreover, recall the definition of $u_n$ in (\ref{fndef}), which was only defined for $n \geq N$. For $0 \leq n \leq N$, let $u_n$ $w_{n}$ be defined implicitly by \eqref{fndef}. The introduction of these objects will not affect the result of our computation.

Observe that
\be
\ds \frac{\gamma^{-1}+K_n(x_0,x_0)}{L_n^2} = \ds \frac{1+\gamma^{-1}}{L_n^2} + R_n  ,
\ee where
\be
\begin{array}{lll}
R_n & = \ds \frac{1}{L_n^2}\ds \sum_{j=1}^n p_j(x_0)^2
&= \ds \sum_{j=1}^{n} \ds \frac{u_{j}^2 (g_{1,j} v_1 + g_{1,j}' r_j v_2)^2}{\lambda_{j+1}^+ \cdots  \lambda_{n}^{+}} \, ,
\end{array}
\ee with the convention that $\lambda_{j+1}^+ \cdots \lambda_{n}^+ = 1$ when $j\geq n$.

Let
\be
\begin{array}{lll}
S_n & =  \ds \frac{ K_{n-1}(x_0,x_0)}{L_{n-1}^2} 
& =  \ds \sum_{j=0}^{n-1} \frac{u_{j}^2 (g_{1,j} v_1 + g_{1,j}' r_j v_2)^2}{\lambda_{j+1}^+ \cdots  \lambda_{n-1}^{+}}.
\label{e12}
\end{array}
\ee

Then
\begin{multline}
\left| \ds \frac{\gamma^{-1}+K_n(x_0,x_0)}{L_n^2} - \ds \frac{\gamma^{-1}+K_{n-1}(x_0,x_0)}{L_{n-1}^2} \right|\\
 \leq  \ds \frac{2(1+\gamma^{-1})}{({\lambda_{n-1}}^+ \cdots  {\lambda_0}^{+})^2} + \left| R_n - S_n \right| .
\label{e12a}
\end{multline}

Recall that $\lambda_n^+ \to \lambda^+ > 1$. Hence,
\be
\ds \sum_n \ds \frac{1}{({\lambda_n}^+ \cdots  {\lambda_0}^{+})^2} = O \left( \ds \sum_{n} \ds \frac{1}{(\lambda^+)^{2}} \right) < \infty .
\ee Thus, the first sum on the right hand side of \eqref{e12a} is summable.

Next, observe that upon rearranging the indices of $S_n$ in \eqref{e12}, we have
\be
\begin{array}{lll}
& |R_n - S_n|  \\
= & \left| \ds \sum_{j=1}^{n}   \ds \frac{u_{j}^2 (g_{1,j} v_1 + g_{1,j}' r_j v_2)^2}{(\lambda_{j+1}^+ \cdots  \lambda_{n}^{+})^2}  - \ds \frac{u_{j-1}^2 (g_{1,j-1} v_1 + g_{1,j-1}' r_{j-1} v_2)^2}{(\lambda_{j}^+ \cdots  \lambda_{n-1}^{+})^2} \right|
\\
\leq & \ds \sum_{j=1}^n \frac{u_{j}^2 |e_{j}-e_{j-1}|}{(\lambda_{j+1}^+ \cdots  \lambda_{n}^{+})^2}   +  \ds \sum_{j=1}^{n} \ds u_{j}^2 e_{j-1} \left| \ds \frac{1}{(\lambda_{j+1}^+ \cdots  \lambda_{n}^{+})^2} - \frac{1}{(\lambda_{j}^+ \cdots  \lambda_{n-1}^{+})^2}          \right| \\
& +
 \ds \sum_{j=1}^n \frac{|u_{j}^2 - u_{j-1}^2|  e_{j-1}}{(\lambda_{j}^+ \cdots  \lambda_{n-1}^{+})^2}  ,
\label{e13} 
\end{array}
\ee where
\be
e_{j}= (g_{2,j} v_1 + g_{2,j}' r_j v_2)^2  .
\label{ejndef}
\ee

Now we proceed to estimate each of the sums on the last line of \eqref{e13}. Recall that $|u_{j} - u_{j-1}| = O(\| G_j - G_{j-1}\|)$. Moreover, $u_{j} \to {u_\infty}$. Hence, $\sum_{j} |u_{j}^2 - u_{j-1}^2| < \infty $. Therefore, for some constant $C$,
\be
\ds \sum_{n=1}^{\infty} \ds \sum_{j=1}^n \frac{|u_{j}^2 - u_{j-1}^2|  e_{j-1}}{(\lambda_{j}^+ \cdots  \lambda_{n-1}^{+})^2} \leq C \left( \ds \sum_{j=1}^{\infty} |u_{j}^2 - u_{j-1}^2| \right) \left(\ds \sum_{j=1}^{\infty} \ds \frac{1}{(\lambda^+)^{2j}}  \right) < \infty .
\ee

Similarly,
\be
\sum_{n=1}^{\infty} \ds \sum_{j=1}^n \frac{u_{j}^2 |e_{j}-e_{j-1}|}{(\lambda_{j+1}^+ \cdots  \lambda_{n}^{+})^2} < \infty 
\ee because
\be \begin{array}{lll}
|e_j - e_{j-1}|  & = O(|r_j - r_{j-1}|) + O(|g_{1,j} - g_{1,j-1}| ) + O(|g_{1,j}' - g_{1,j-1}'| ) \\
&=  O(\| G_{j}-G_{j-1}\| ) .
\end{array}
\ee

Finally, we consider the second sum on the right hand side of \eqref{e13}. Note that
\be
|{\lambda_{j}^+}^2 - {\lambda_n^+}^2| \leq C |{\lambda_j^+} - {\lambda_n^+} | \leq C \ds \sum_{k=j}^{n-1} |\lambda_{k+1}^+ - \lambda_k^+| \, ,
\ee where $C$ is a positive constant independent of $j$ and $n$. As a result,
\be \begin{array}{lll}
 & \ds \sum_{n=1}^{\infty} \ds \sum_{j=1}^{n} \left| \ds \frac{1}{(\lambda_{j}^+ \dots \lambda_{n-1}^+)^2} - \ds \frac{1}{(\lambda_{j+1}^+ \dots \lambda_{n}^+)^2}  \right|   \\
  =&  \ds \sum_{n=1}^{\infty} \ds \sum_{j=1}^{n} \ds \frac{|{\lambda_{j}^+}^2 - {\lambda_{n}^+}^2|}{(\lambda_{j}^+ \dots \lambda_{n}^+)^2}  \\
 \leq & C \ds \sum_{n=1}^{\infty} \ds \sum_{j=1}^{n} \frac{1 }{(\lambda_{j+1}^+ \dots \lambda_{n}^+)^2}\ds \sum_{k=j}^{n-1} |\lambda_{k+1}^+ - \lambda_k^+|.
%
%
\end{array} 
\label{e14}
\ee 

Consider a fixed $k \in \mathbb{N}$. We count the coefficients of $|\lambda_{k+1}^+ - \lambda_k^+|$ in the last line of \eqref{e14}. From the equation, we know that $j \leq k <n$. Thus, the coefficient is
\be \begin{array}{lll}
& \ds \sum_{n=k+1}^{\infty} \ds \sum_{j=1}^{k} \frac{1 }{(\lambda_{j+1}^+ \dots \lambda_{n}^+)^2} 
\\ 
= & \ds \sum_{j=1}^k \left( \ds \sum_{n=k+1}^\infty \frac{1}{(\lambda_{j+1}^+ \cdots \lambda_{n}^+)^2 } \right) \\
= & \left( \ds \sum_{j=1}^{k} \ds \frac{1}{(\lambda_2^+ \cdots \lambda_{j}^+)^2} \right) \left( \ds \sum_{n=k+1}^\infty \frac{1}{(\lambda_{k+1}^+ \cdots \lambda_n^+)^2 } \right),
\end{array}
\ee which is bounded above by a constant $B$ independent of $k$. Therefore,
\be
\eqref{e14} \leq C B \ds \sum_{k=1}^{\infty} |\lambda_{k+1}^+ - \lambda_k^+| < \infty .
\ee Going back to \eqref{e13}, we conclude that $\sum_n |R_n-S_n| < \infty$. Hence, $\tilde{a}_n$ is of bounded variation.

Finally, we are going to show that $\tilde{b}_n$ is of bounded variation.
Since $a_n$ is of bounded variation, by \eqref{tbn}, it suffices to show that $(h_n)_n$  is of bounded variation, $h_n$ as defined in \eqref{hndef}.

By \eqref{pnasymp},
\be
h_n =\ds \frac{ a_{n+1} \lambda_{n+1}^+ L_n^2 u_{n+1} u_n }{\gamma^{-1} + K_n(x_0, x_0)} \left( g_{1,n+1} v_1 + g_{1,n+1}' r_{n+1} v_2  \right) \left( g_{1,n} v_1 + g_{1,n}' r_n v_2  \right) .
\label{e100}
\ee

Thus, it boils down to proving that $L_n^2/(\gamma^{-1} + K_n(x_0,x_0))$ is of bounded variation. To that end, we use the same argument following \eqref{taundef}.

\section{Acknowledgements}
It is my pleasure to thank Professor Jeff Geronimo and Professor Doron Lubinsky for very helpful discussions.


\begin{thebibliography}{8}

\bibitem{borg} G. Borg, \emph{Eine Umkehrung der Sturm-Liouvilleschen Eigenwertaufgabe. Bestimmung der Differentialgleichung durch die Eigenwerte},
Acta Math. \textbf{78} (1946), 1--96. 

\bibitem{cm1} A. Cachafeiro and F. Marcell\'an, \emph{Orthogonal polynomials and jump modifications}, in ''Orthogonal Polynomials and Their Applications'', (Segovia, 1986), pp 236--240, Lecture Notes in Math., 1329, Springer, Berlin, 1988.

\bibitem{cesaro} E. Ces\`aro and O. Stolz, \emph{http://en.wikipedia.org/wiki/Stolz-Ces\'aro\_theorem}.

\bibitem{chihara} T. Chihara, \emph{An Introduction to Orthogonal Polynomials}, Mathematics and its applications \textbf{13}, Gordon and Breach, New York, 1978.


\bibitem{dks} D. Damanik, R. Killip, B. Simon, \emph{Perturbations of orthogonal polynomials with periodic recursion coefficients}, Annals of Math., \textbf{171} (2010), No. 3, 1931--2010



\bibitem{dombrowski} J. Dombrowski, \emph{Quasitriangular Matrices}, Proc. Amer. Math. Soc. \textbf{69} (1978), 95--96.

\bibitem{gelfand} I. M. Gel'fand and B. M. Levitan, \emph{On the determination of a differential equation from its spectral function}, Amer. Math. Soc. Transl. (2) \textbf{1} (1955), 253--304; Russian original in Izvestiya Akad. Nauk SSSR. Ser. Mat. \textbf{15} (1951), 309--360.

\bibitem{birthday} F. Gesztesy, P. Deift, C. Galvez, P. Perry \& W. Schlag, editors: \emph{Spectral Theory and Mathematical Physics: A Festschrift in Honor of Barry Simon's 60th Birthday}, Proc. Sympos. Pure Math. 76, Amer. Math. Soc., Providence, RI, 2007.

\bibitem{gv} J. S. Geronimo and W. Van Assche, \emph{Orthogonal polynomials with asymptotically periodic recurrence coefficients}, J. Approx. Theory \textbf{46} (1986), 251--283
\bibitem{golinskiinevai} L. Golinskii and P. Nevai, \emph{Szeg\H o difference equations, transfer matrices and orthogonal polynomials on the unit circle}, Comm. Math. Phys. \textbf{223} (2001), 223--259.

\bibitem{geronimus1} Ya. L. Geronimus, \emph{Polynomials Orthogonal on a Circle and Their Applications}, Amer. Math. Soc. Translation \textbf{104} (1954) 79pp.

\bibitem{geronimus2}  Ya. L. Geronimus, \emph{Orthogonal Polynomials: Estimates, Asymptotic Formulas, and Series of Polynomials Orthogonal on the Unit Circle and on an Interval}, Consultants Bureau, New York, 1961.

\bibitem{kooman1} R. J. Kooman, \emph{Asymptotics behaviour of solutions of linear recurrences and sequences of M\"obius-transformations}, J. Approx. Theory \textbf{93} (1998), 1--58.

\bibitem{kooman2} R. J. Kooman, \emph{Decomposition of matrix sequences}, Indag. Math. (N. S.) \textbf{5} (1944) 61--79.

\bibitem{lubinsky1} D. S. Lubinsky, \emph{A survey of general orthogonal polynomials for weights on finite and infinite intervals}, Acta Appl. Math. \textbf{10} (1987), 237--296. 


\bibitem{mnt} A. M\'at\'e, P. Nevai, V. Totik, \emph{Asymptotics for orthogonal polynomials defined by a recurrence relation}, Constr. Approx. 1 \textbf{3} (1985), 231--248.

\bibitem{nevai1} P. Nevai, \emph{Orthogonal Polynomials}, Mem. Amer. Math. Soc. \textbf{18} (1979), no. 213, 185 pp.

\bibitem{nevai2} P. Nevai, \emph{Geza Freud orthogonal polynomials and Christoffel functions. A case study}, J. Approx. Theory \textbf{48} (1986), 3--167. 


\bibitem{ntz} P. Nevai, V. Totik, J. Zhang, \emph{Orthogonal polynomials: their growth relative to their sums}, J. Approx. Theory \textbf{67} (1991), 215--234.


\bibitem{simon1} B. Simon, \emph{Orthogonal polynomials on the unit circle. Part 1: Classical Theory}, AMS Colloquium Series, Amer. Math. Soc., Providence, RI, 2005.

\bibitem{simon2} B. Simon, \emph{Orthogonal polynomials on the unit circle. Part 2: Spectral Theory}, AMS Colloquium Series, Amer. Math. Soc., Providence, RI, 2005.

\bibitem{simon3} B. Simon, \emph{The Christoffel-Darboux kernel}, Proc. Sympos. Pure Math., \textbf{79}, Amer. Math. Soc., Providence, RI, 2008.


\bibitem{totik} W. Totik, \emph{Orthogonal polynomials}, Surveys in Approx. Theory \textbf{1} (2005), 70--125.

\bibitem{wigner} J. von Neumann and E. Wigner, \emph{\"Uber merkw\"urdige diskrete Eigenwerte}, Phys. Z. \textbf{30} (1929), 465--467.


\bibitem{szwarc} R. Szwarc, \emph{Chain sequences and compact perturbations of orthogonal polynomials}, Math. Zeits. \textbf{217} (1994), 57--71.

\bibitem{wong1} M.-W. L. Wong, \emph{Generalized bounded variation and inserting point masses}, Const. Approx. \textbf{30} (2009), 1--15.


\bibitem{wong2} M. -W. L. Wong, \emph{A formula for inserting point masses}, J. Comput. and Appl. Math., \textbf{233} (2009), 852--885.


\bibitem{wong3} M. -W. L. Wong, \emph{Asymptotics of orthogonal polynomials and point perturbation on the unit circle}, J. Approx. Theory \textbf{162} (2010), 1294--1321.

\bibitem{wong4} M. -W. L. Wong, \emph{Point mass insertion on the real line and non-exponential decay of the recurrence coefficients}. To appear in Math. Nachr.



\end{thebibliography}
\end{document}